\documentclass[11pt]{article}
\usepackage{amssymb}
\usepackage{amsmath,epsf,wrapfig,multicol,epic,eepic,trig,mathrsfs,dsfont,color}
\usepackage[dvipdfm]{graphicx}
\newtheorem{thm}{Theorem}
\newtheorem{cor}[thm]{Corollary}
\newtheorem{lem}[thm]{Lemma}

\newtheorem{conj}[thm]{Conjecture}
\newtheorem{definition}[thm]{Definition}
\newtheorem{rmk}[thm]{Remark}
\newtheorem{example}[thm]{Example}
\makeatletter
\newenvironment{pf}{{\bf Proof.} }
\makeatother
\makeatletter
\newenvironment{pot}{{\bf Proof of Theorem \ref{thm1}. } }
\makeatother

\newcommand{\Z}{\mathbb{Z}}

\newcommand{\qed}{\hbox{\rule[-2pt]{3pt}{6pt}}}
\def\dl{\langle\mspace{-5mu}\langle}
\def\dr{\rangle\mspace{-5mu}\rangle}

\input{naoko15Fg.tex}

\begin{document}
\title
{Surface pole bracket polynomials of virtual knots and twisted knots}
\author{Naoko Kamada\\
Graduate School of Natural Sciences, Nagoya City University, \\Mizuho-ku, Nagoya, Aichi 467-8501, Japan}

\maketitle

\begin{abstract} 
Dye and Kauffman defined surface bracket polynomials for virtual links by use of surface states, and  
found a relationship between the surface states and the minimal genus of a surface in which a virtual link diagram is realized.  They and Miyazawa independently defined  
a multivariable polynomial invariant of virtual links.   This invariant is deeply related to the surface states.  
In this paper,  we introduce the notion of surface pole bracket polynomials for link diagrams in closed surfaces, as a  generalization of surface bracket polynomials by Dye and Kauffman.  The polynomials induce the invariant of twisted links defined by the author before as a generalization of Dye, Kauffman and Miyazawa's polynomial invariant. 
Furthermore we discuss a relationship between curves in surface pole states and variables of the polynomial invariant.
\end{abstract}

\section{Introduction}
Virtual knot theory is a generalization of knot theory which is based on  Gauss chord diagrams and link diagrams on closed oriented surfaces  \cite{rkauD}. Virtual links correspond to  stable equivalence classes of  links in oriented 3-manifolds which are line bundles over closed  oriented surfaces (cf. \cite{rCKS, rkk}). A twisted link defined by Bourgoin \cite{rBor} is  an extension of the notion of vital links.  Twisted links correspond to stable equivalence classes of links in oriented 3-manifolds which are line bundle over closed surfaces which are possibly non-orientable surfaces \cite{rBor}. 

\begin{figure}[h]
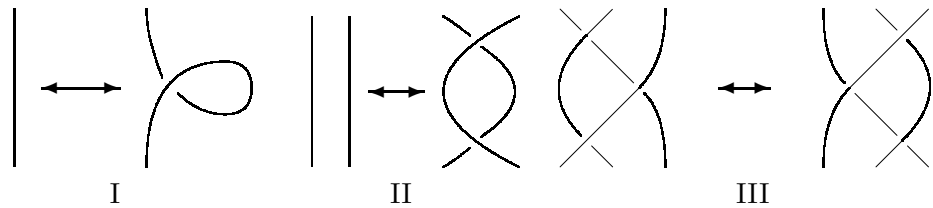

\begin{center}
\begin{tabular}{ccc}
\rmoveio{.7mm}{1}&\rmoveiio{.5mm}{1}&\rmoveiiio{.7mm}{1}\\
I&II&III\\
\end{tabular}
\caption{Reisdemeister moves}\label{fig:movesR}
\end{center}
\end{figure}

\begin{figure}[h]
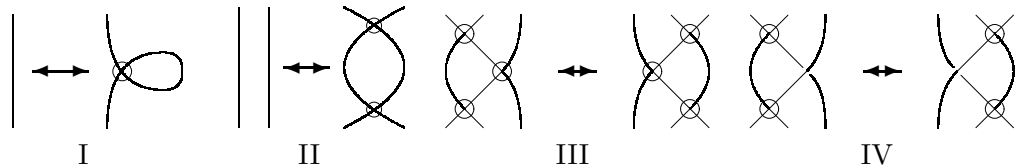

\begin{center}
\begin{tabular}{cccc}
\rmovevio{.5mm}{1}&\rmoveviio{.4mm}{1}&\rmoveviiio{.5mm}{1}&\rmovevivo{.5mm}{1}\\
I&II&III&IV\\
\end{tabular}
\caption{Virtual Reidemeister moves}\label{fig:movesV}
\end{center}
\end{figure}

A {\it virtual link diagram\/} is a link diagram which may have {\it virtual crossings\/}, which are encircled crossings without over-under information.  
A {\it virtual link\/}  is an equivalence class of a virtual link diagram by Reidemeister moves and virtual Reidemeister moves depicted in Figures~\ref{fig:movesR} and \ref{fig:movesV}. We call these moves {\it generalized Reisemeister moves\/}. 

\begin{figure}[h]
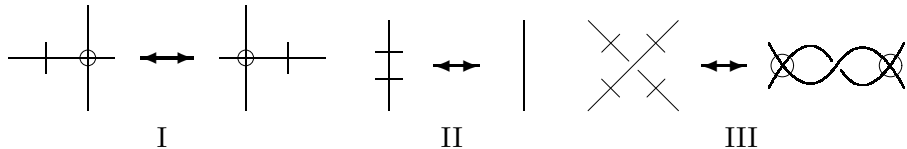

\begin{center}
\begin{tabular}{ccc}
\rmoveti{.7mm}&\rmovetiio{.6mm}&\rmovetiiio{.6mm}{1}\\
I&II&III\\
\end{tabular}
\caption{Twisted Reidemeister moves}\label{fig:movesT}
\end{center}
\end{figure}

A {\it twisted link diagram\/} is a virtual link diagram which may have {\it bars\/} on arcs. 
A {\it twisted link\/}  is an equivalence class of a twisted link diagram by Reidemeister moves, virtual Reidemeister moves and twisted Reidemeister moves in Figures~\ref{fig:movesR}, \ref{fig:movesV} and \ref{fig:movesT}. We call these moves  {\it extended Reisemeister moves\/}. 

\vspace{.2cm}

Bourgoin introduced  the Jones polynomials ($f$-polynomials) for twisted links and a group invariant called the twisted knot group \cite{rBor}. The author introduced a twisted  quandle for twisted links \cite{rkamF}. 
For a twisted link $L$, it is an interesting and important problem to determine an irreducible representative or to determine the minimum genus of a surface $F$ in which a diagram of $L$ is realized.  
Surface bracket polynomials of virtual links are defined by Dye and Kauffman \cite{rDKB} by use of surface states which are obtained from a link diagram in a closed oriented surface in which a diagram of $L$ is realized. 

The following conjecture is due to Kauffman and  Przytycki.  
\begin{conj}
For a virtual knot $L$, if a diagram of $L$ is realized in a surface of the minimal genus, then this fact is detected by the surfce bracket polynomial.
\end{conj}

H. A. Dye and L. H. Kauffman \cite{rDKA}, and Y. Miyazawa \cite{rMiyaB} independently, defined 
a multivariable polynomial invariant of virtual links, which we call the 
DKM polynomial.   Dye and Kauffman showed that this invariant is deeply related to the surface states for link diagrams on closed oriented surfaces.  
In this paper,  we introduce the notion of surface pole bracket polynomials for link diagrams in closed surfaces, as a  generalization of surface bracket polynomials by Dye and Kauffman.  The polynomials induce the invariant of twisted links defined by the author in \cite{rkamF}  as a generalization of the DKM polynomial invariant. 
Then we discuss a relationship between curves in surface pole states and variables of the polynomial invariant.

\section{Link diagram realizations of twisted links}

\begin{figure}[h]
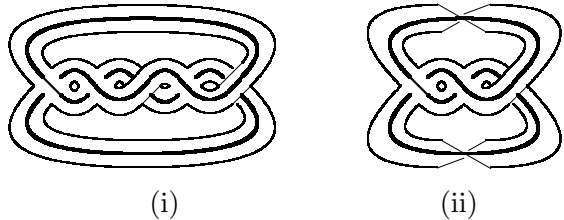

\begin{center}
\begin{tabular}{cc}
\absexcc{.6mm}
&
\abstwtexcc{.6mm}\\
(i)&(ii)\\
\end{tabular}
\caption{Examples of abstract link diagrams}\label{fg:exabstdiag}
\end{center}
\end{figure}

An {\it abstract link diagram}  is a pair $(\Sigma, D_{\Sigma})$ 
of a compact, possibly non-orientable surface $\Sigma$ and a link diagram $D_{\Sigma}$ in $\Sigma$ such that $|D_{\Sigma}|$ is a deformation retract of $\Sigma$, where $|D_{\Sigma}|$ is the subset of $\Sigma$ obtained from $D_{\Sigma}$ by replacing each crossing with a $4$-valent vertex.  
Two examples of abstract links are depicted in Figure~\ref{fg:exabstdiag}.  The surface $\Sigma$ in (ii) of the figure is a non-orientable surface. 

Let $(\Sigma_1, D_{\Sigma_1})$ and  $(\Sigma_2, D_{\Sigma_2})$ be abstract link diagrams, where $\Sigma_i$ is a compact surface and $D_{\Sigma_i}$ is a link diagram in $\Sigma_i$ for each $i \in \{1,2\}$. 
If there are embeddings $f_1: \Sigma_1 \to F$ and $f_2: \Sigma_2 \to F$ to a common closed surface $F$ such that 
 $f_2(D_{\Sigma_2})$ is obtained from $f_1(D_{\Sigma_1})$ by a Reidemeister move in $F$, then we say that $(\Sigma_2, D_{\Sigma_2})$ is obtained from $(\Sigma_1, D_{\Sigma_1})$ by an {\it abstract Reidemeister move}.  
 Two abstract link diagrams are said to be {\it equivalent} if they are related by a finite sequence of abstract Reidemeister moves.  
We call  equivalence classes of abstract link diagrams  {\it  abstract links}. 
(This notion was introduced in \cite{rkk} for the case where surfaces $\Sigma$'s and $F$'s are oriented.  Note that, in this paper, we do not assume that these surfaces are orientable.) 

\begin{thm}[Bourgoin \cite{rBor}] 
There is a map from the family of twisted link diagrams to that of abstract link diagrams such that it induces a bijection from the family of twisted links to that of abstract links.
\end{thm}

\begin{figure}[h]
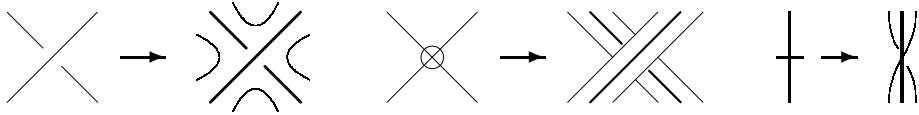

\begin{center}
\begin{tabular}{ccc}
\abscrs{.6mm}&\absvcrs{.6mm}&\abstwt{.6mm}
\end{tabular}
\caption{Twisted link diagram and abstract link diagram}\label{fg:twstvsabst}
\end{center}
\end{figure}

\begin{figure}[h]
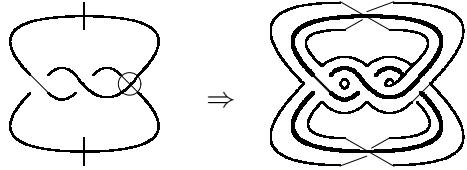

\begin{center}
\begin{tabular}{ccc}
\twtexcc{.6mm}&\raisebox{.9cm}{$\Rightarrow$}&
\abstwtexcc{.6mm}
\end{tabular}
\caption{Abstract link diagram associated  with twisted link diagram}\label{fg:texwstvsabst}
\end{center}
\end{figure}

This map is depicted in Figure~\ref{fg:twstvsabst}. We call the abstract link diagram obtained this way the {\it abstract link diagram associated with $D$}. For example, see Figure~\ref{fg:texwstvsabst}. 

\begin{figure}[h]
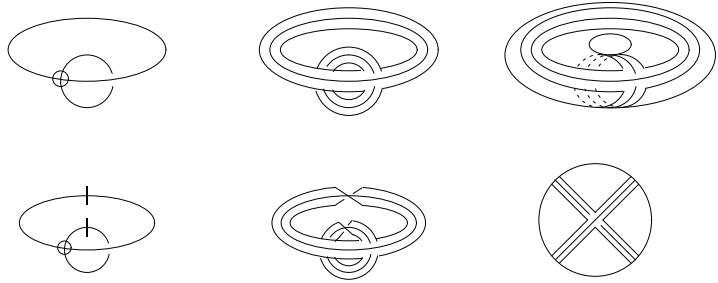

\begin{center}
\begin{tabular}{ccc}
\vir{.7mm}
&
\absvirb{.7mm}
&
\survirb{.7mm}
\\
\twstb{.6mm}
&
\abstwstb{.6mm}
&
\surtwstb{.5mm}
\end{tabular}
\caption{Examples of link diagram realizations}\label{fg:exlinkdiagreal}
\end{center}
\end{figure}

A pair $(F,D_F)$ of a closed surface $F$ and a link diagram $D_F$ in $F$ is called a {\it link diagram realization} of a twisted link diagram $D$ if there is an embedding $f: \Sigma \to F$ such that $f(D_{\Sigma})=D_F$, where $(\Sigma, D_{\Sigma})$ is the abstract link diagram associated with $D$. 
For example, see Figure~\ref{fg:exlinkdiagreal}, where the link diagram realization depicted in the bottom right is a link diagram in a projective plane.

\section{Surface pole states}

\begin{figure}[h]
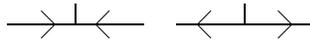

\begin{center}
\pole{.9mm}
\caption{A sink pole (I-pole) and a source pole (O-pole)}\label{fig:pole}
\end{center}
\end{figure}

A  {\it pole curve} in a surface $F$ is a simple closed curve with (or without) poles in $F$ as in Figure~\ref{fig:pole}, where we call the pole on the left side of the figure a {\it sink pole} (or an {\it I-pole}) and the one on the right side a {\it source pole} (or an {\it O-pole}). 
A collection of mutually disjoint pole curves in a surface $F$ is called a {\it pole curve link} in $F$.  

\begin{figure}[h]
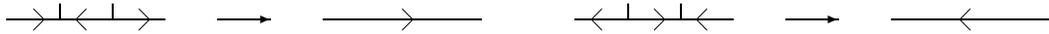

\vspace{.3cm}
\begin{center}
\begin{tabular}{cc}
\polereductA{.7mm}&\polereductB{.7mm}
\end{tabular}
\caption{Pole reductions}\label{fg:polereduct}
\end{center}
\end{figure}

The local moves depicted in Figure~\ref{fg:polereduct} are called {\it pole reductions}.  

\begin{figure}[h]
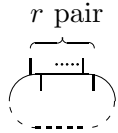

\vspace{.3cm}
\begin{center}
\polecntA{.7mm}
\caption{An irreducible pole curve}\label{fg:irrpolediag}
\end{center}
\end{figure}

A pole curve is said to be {\it irreducible} if one cannot apply any pole reduction to it. 
An irreducible pole curve looks like 
as in Figure~\ref{fg:irrpolediag}, where there is no pole in the broken line of the curve. 
A pole curve link is said to be {\it irreducible} if every component is irreducible. 

For a pole curve $l$ in $F$, applying pole reductions, we obtain an irreducible pole curve in $F$.  We denote it by 
$\bar{l}$.  The {\it index} of $l$ is the half of the numbers of poles of the irreducible pole curve $\bar{l}$ obtained from $l$, which is denoted by $\iota (l)$.

\begin{figure}[h]
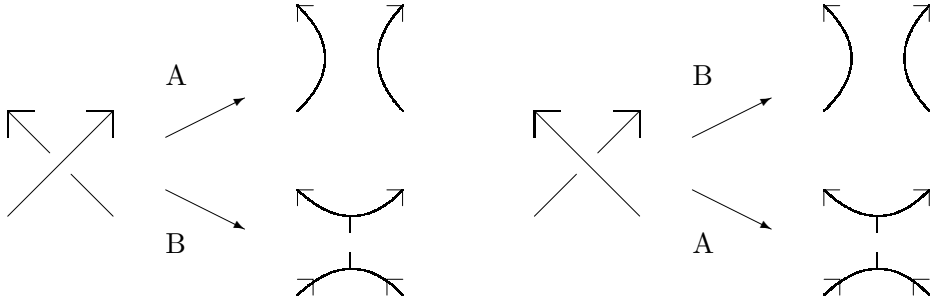

\vspace{.3cm}
\begin{center}
\splicepb{.7mm}
\caption{A-splices and B-splices}\label{fg:splice}
\end{center}
\end{figure}

Let $D$ be a twisted link diagram and $(F,D_F)$  a link diagram realization of $D$ in a surface $F$. 
The local replacement at a crossing of a link diagram in $F$ illustrated in Figure~\ref{fg:splice} 
is called an {\it A-splice} or a {\it B-splice}.  

Applying A- or B-splice at each crossing of $D_F$, we have a pole curve link in $F$, say $s$.  
The pair $(F, s)$ or $s$ is called a 
 {\it surface pole state} of $(F,D_F)$.  
Moreover, applying pole reductions of all pole curves in $s$, we have an irreducible pole curve link in $F$.  We denote it by $\bar{s}$, and call the pair $(F, \bar{s})$ an {\it irreducible pole state} associated with $(F, s)$.  
(An irreducible pole curve link $\bar{s}$ obtained from $s$ is determined uniquely up to equivalence in the sense of Definition~\ref{defn:equivpolelink} (Remark~\ref{remark:equivpolelink}).)

When $D$ is a virtual knot and $(F, D_F)$ is a link diagram realization of $D$ in a closed orinted surface $F$, the notions  of pole diagrams and surface pole states are essentially due to Dye and Kauffman \cite{rDKA}.

\begin{thm}[Dye and Kauffman \cite{rDKA}]\label{thm2}
Let $D$ be a virtual link diagram and $(F, D_F)$  a link diagram realization of $D$ in a closed orinted surface $F$.  
Let $l$ be a pole curve of a surface pole state $(F, s)$ of $(F, D_F)$.  
If the index $\iota(l)$ is positive, then $l$ is an essential curve in $F$.
\end{thm}

Theorem~\ref{thm2} holds for twisted links.

\begin{thm}\label{thm1}
Let $D$ be a twisted link diagram and $(F, D_F)$  a link diagram realization of $D$. 
Let $l$ be a pole curve of a surface pole state $(F, s)$ of $(F, D_F)$.  
If the index $\iota(l)$ is positive, then $l$ is not a separating curve in $F$.  
\end{thm}

Note that if a simple closed curve in $F$ is not a separating curve in $F$ then it is an essential curve in $F$.  
From this theorem, we have the following.

\begin{cor}\label{cor1}
Let $D$ be a virtual link diagram and $(F, D_F)$  a link diagram realization of $D$.    
Let $l$ be a pole curve of a surface pole state $(F, s)$ of $(F, D_F)$.  
If the index $\iota(l)$ is positive, then $l$ is not a separating curve in $F$.
\end{cor}

\section{Proof of Theorem~\ref{thm1}}

\begin{lem}\label{lem1}
For a pole curve $l$ in $F$, the number of I-poles on $l$ is equal to that of O-poles. 
\end{lem}

\begin{pf}
Since I-poles and O-poles appear adjacently on $l$, the number of I-poles on $l$ is equal to that of O-poles.  \qed
\end{pf}

\begin{lem}\label{lem2}
Let $D$ be a twisted link diagram and $(F, D_F)$  a link diagram realization of $D$. 
For a suface pole state $(F, s)$ of $(F, D_F)$, let  $R$  be a connected component of $F\setminus s$. The number of I-poles on $s$ in $R$ is equal to that of O-poles.
\end{lem}

\begin{pf}
By a splice yielding a pair of poles, two poles come up in the same component of $F\setminus s$. 
Since one is an I-pole and the other is an O-pole in such a pair, 
the number of I-poles in $R$ is equal to that of O-poles.\qed 
\end{pf}

By a pole reduction, a pair of an I-pole and an O-pole is reduced in the same component of $F\setminus s$.  Therefore we have the following from the above lemma. 

\begin{lem}\label{lem3}\quad 
Let $D$ be a twisted link diagram and $(F, D_F)$ a link diagram realization of $D$. 
Let $(F, \bar{s})$ be an 
irreducible surface pole state obtained from a surface pole state $(F, s)$ of $(F, D_F)$, and let  $R$  be a connected component of $F\setminus \bar{s}$ $( = F \setminus s)$. The number of I-poles on $\bar{s}$ in $R$ is equal to that of O-poles.
\end{lem}

\begin{figure}[h]
\vspace{.3cm}
\begin{center}
\includegraphics[width=8cm]{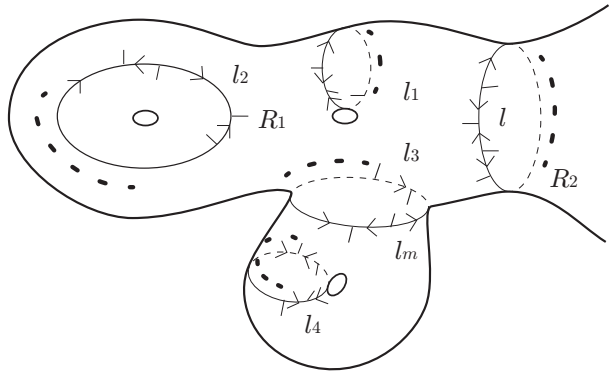}
\caption{Irreducible surface pole state}\label{fg:statesurf}
\end{center}
\end{figure}

\begin{figure}[h]
\vspace{.3cm}
\begin{center}
\includegraphics[width=5cm]{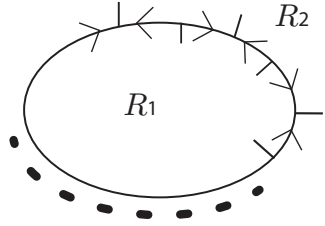}
\caption{Irreducible pole curve whose index is positive}\label{fg:sepcurve}
\end{center}
\end{figure}

\begin{pot}
Let $(F, \bar{s})$ be an irreducible surface pole state of a surface pole state $(F, s)$.  
Let $\bar{l}$ be the irreducible pole curve of $(F, \bar{s})$ obtained from a pole curve $l$ of $s$ such that $\iota(l)$ is positive.   Assume that $F\setminus \bar{l}$ consists of two components, say $R_1$ and $R_2$, and we show a contradiction.   
Let $l_1,\dots,l_m$ be irreducible pole curves of $(F, \bar{s})$ in $R_1$ and $n_1(I), \dots , n_m(I)$ (or $n_1(O), \dots, n_m(O)$, resp.) be the numbers of I-poles (or O-poles, resp.) on $l_1, \dots, l_m$ (see Figure~\ref{fg:statesurf}). 
Let $n^i(I)$ (or $n^i(O)$, resp.) be the number of I-poles (or O-poles, resp.) on $\bar{l}$ in $R_i$ for $i=1,2$. 
By Lemma~\ref{lem1}, we have $n_k(I)=n_k(O)$ for $k=1,\dots, m$ and $n^1(I)+n^2(I)=n^1(O)+n^2(O)$.  
On the other hand, 
considering all connected components of $F \setminus \bar{s}$ contained in $R_1$, 
we see by Lemma~\ref{lem3} that 
$$\displaystyle \sum_{k=1}^m n_k(I)+n^1(I)=\sum_{k=1}^m n_k(O)+n^1(O).$$ 
Thus we have that $n^1(I)=n^1(O)$ and $n^2(I)=n^2(O)$. 
Since $\iota(l)$ is positive, all O-poles (or all I-poles, resp.) on $\bar{l}$ are in $R_1$ and all I-poles (or all O-poles, resp.) on $\bar{l}$ are in $R_2$ (cf. Figure~\ref{fg:sepcurve}). This implies that $n^1(O)$ is positive and $n^1(I)=0$ (or $n^1(I)$ is positive and $n^1(O)=0$, resp.). It contradicts $n^1(I)=n^1(O)$. \qed
\end{pot}

\section{Surface pole bracket polynomials}
\label{sect:surfacepolebracketpoly}

In this section we introduce the notion of a surface pole bracket polynomial.  

Let $F$ be a closed, possibly non-orientable, surface.  

\begin{definition}\label{defn:equivpolelink}{\rm 
Two pole links $s$ and $s'$ in $F$ are {\it equivalent} if $s'$ is obtained from $s$ by an ambient isotopy of $F$ and by changing the orientations of some (or none) of curves without poles.  
}\end{definition}

\begin{rmk}\label{remark:equivpolelink}{\rm 
Let $l$ be a pole curve in $F$ with some poles.  If $l$ has some poles and $\iota(l)=0$, then by pole reductions we obtain a pole curve without poles, i.e., a simple loop in $F$.  It is an irreducible pole curve $\overline{l}$ obtained from $l$.  
The orientation of $\overline{l}$ is not determined uniquely from $l$.  Actually, changing how to apply pole reductions, we may obtain a simple loop with the opposite orientation. Thus, when we consider  $\overline{l}$ for a pole curve $l$ with $\iota(l)=0$, it is more natural to ignore the orientation of $\overline{l}$.   By introducing Definition~\ref{defn:equivpolelink}, 
we may say that 
for a pole curve $l$ in $F$,  an irreducible pole curve $\overline{l}$ obtained from $l$ by pole reductions is uniquely determined up to equivalence, and that for a pole curve link $s$ in $F$,  an irreducible pole curve link $\overline{s}$ obtained from $s$ is uniquely determined up to equivalence. 
}\end{rmk}

Let $P(F)$ be the family of all equivalence classes of pole curve links in $F$.  (We allow the empty set as a pole curve link in $F$.) Let $\Z[A, A^{-1}] P(F)$ be the free module generated by $P(F)$ over the Laurent polynomial ring $\Z[A, A^{-1}]$. 

For a pole curve link $s$ in $F$, we denote by ${\rm ess}(s)$ (or by ${\rm iness}(s)$, resp.) the subset of $s$ consisting of pole curves that are essential loops in $F$ (or inessential loops in $F$, resp.).  Note that $s = {\rm ess}(s) \amalg 
{\rm iness}(s)$.  By $\overline{s}$ we denote an irreducible pole curve link obtained from $s$ by pole reductions.

\begin{definition}{\rm 
For a pair $(F, D_F)$ of a closed surface $F$ and a link diagram $D_F$ in $F$, 
the {\it surface pole bracket polynomial} $\langle (F,D_F) \rangle$ of $(F, D_F)$ is define by 

$$\langle (F,D_F) \rangle=\sum_{s}A^{\natural (s)}(-A^2-A^{-2})^{\sharp ({\rm iness}(s))}[\overline{{\rm ess}(s)}] 
\quad \in \Z[A, A^{-1}] P(F),$$
where in the summation $s$ runs all over surface pole states $(F, s)$ of $(F, D_F)$, 
$\natural (s)$ is the number of A-splices minus that of B-splices obtaining the state $(F, s)$, 
$\sharp ({\rm iness}(s))$ is the number of components of ${\rm iness}(s)$, 
and $[\overline{{\rm ess}(s)}]$ is the equivalence class of an irreducible pole curve link $\overline{{\rm ess}(s)}$ 
obtained from ${\rm ess}(s)$ by pole reductions.  
}\end{definition}

When $F$ is orientable, this notion is essentially introduced to Dye and Kauffman \cite{rDKA}. 
The surface pole bracket polynomial is not an invariant of a link in the thickened surface.  However it might be useful for study such links,  virtual links and twisted links as discussed in \cite{rDKA}.

\begin{example}{\rm 
Let $F$ be a torus and let $D_F$ be a diagram in $F$ illustrated in the top right of Figure~\ref{fg:exsufstate}. 
It has $3$ crossings and there are $8$ surface pole states as depicted in the figure.  In the figure, for each state $s$, $A^{\natural (s)}$ is indicated. 

\begin{figure}[h]
\begin{center}
\includegraphics[width=12cm]{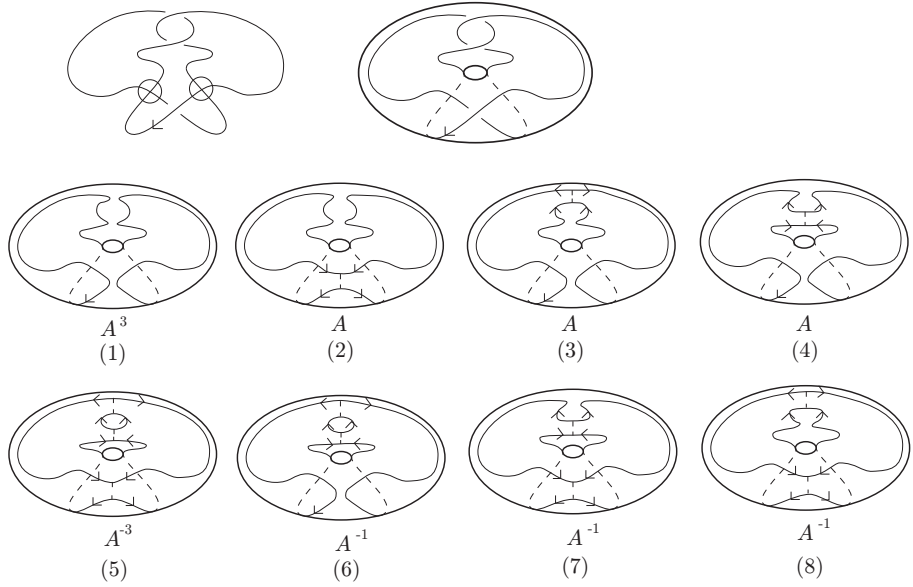}
\caption{Example of a surface pole bracket polynomial}\label{fg:exsufstate}
\end{center}
\end{figure}

Let $s_1, \dots, s_8$ be the states depicted in (1)--(8) of Figure~\ref{fg:exsufstate}  and let 
$L_0, \dots, L_3$ be pole curve links depicted in (i)--(iv) of Figure~\ref{fg:exsufstate4}, respectively. 
Then we have 
$$
\begin{array}{lll}
{\natural (s_1)} = 3,  &   {\sharp ({\rm iness}(s_1))}= 0, & [\overline{{\rm ess}(s_1)}]=[L_1],  \\
{\natural (s_2)} = 1,  &   {\sharp ({\rm iness}(s_2))}= 1, & [\overline{{\rm ess}(s_2)}]=[L_0],  \\
{\natural (s_3)} = 1,  &   {\sharp ({\rm iness}(s_3))}= 1, & [\overline{{\rm ess}(s_3)}]=[L_0] ,  \\
{\natural (s_4)} = 1,  &   {\sharp ({\rm iness}(s_4))}= 1, & [\overline{{\rm ess}(s_4)}]=[L_0] ,  \\
{\natural (s_5)} = -3,  &   {\sharp ({\rm iness}(s_5))}= 1, & [\overline{{\rm ess}(s_5)}]=[L_2],  \\
{\natural (s_6)} = -1,  &   {\sharp ({\rm iness}(s_6))}= 2, & [\overline{{\rm ess}(s_6)}]=[L_0],  \\
{\natural (s_7)} = -1,  &   {\sharp ({\rm iness}(s_7))}= 0, & [\overline{{\rm ess}(s_7)}]=[L_2],  \\
{\natural (s_8)} = -1,  &   {\sharp ({\rm iness}(s_8))}= 0, & [\overline{{\rm ess}(s_8)}]=[L_2].  \\
\end{array}
$$
Thus $\langle (F,D_F) \rangle= (3A(-A^2-A^{-2}) + A^{-1}(-A^2-A^{-2})^2) [L_0] + A^3[L_1] + A^{-3}(-A^2-A^{-2})[L_2]$ 
$= (2A-A^{-3})(-A^2-A^{-2}) [L_0] + A^3[L_1] + A^{-3}(-A^2-A^{-2})[L_2]$. 

\begin{figure}[h]
\begin{center}
\includegraphics[width=12cm]{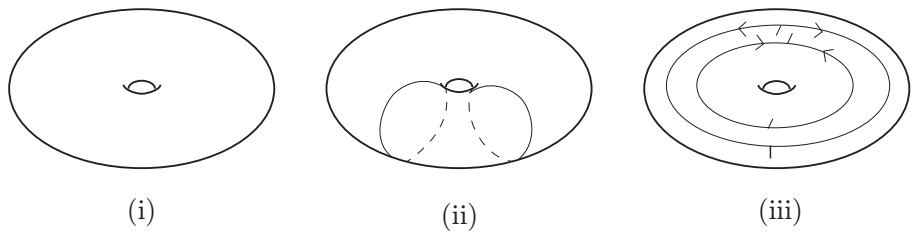}
\caption{Pole curve links}\label{fg:exsufstate4}
\end{center}
\end{figure}

}\end{example}

\begin{definition}{\rm 
Let $(F, D_F)$ be a pair of a closed surface $F$ and a link diagram $D_F$ in $F$.  
The {\it double bracket polynomial} $\dl (F, D_F) \dr$ of $(F, D_F)$ is defined by 
$$\begin{array}{ll}
\dl (F, D_F) \dr=\sum_{s}A^{\natural (s)}(-A^2-A^{-2})^{\sharp ({\rm iness}(s))}&M^{\sharp({\rm non\mbox{-}ori}(s))}\prod_{l \in s}d_{\iota(l)} \\
&\quad \in \Z[A, A^{-1}, M, d_1, d_2, \dots],
\end{array}$$
where in the summation $s$ runs all over surface pole states $(F, s)$ of $(F, D_F)$, 
$\natural (s)$ is the number of A-splices minus that of B-splices obtaining the state $(F, s)$, 
$\sharp ({\rm iness}(s))$ is the number of components of ${\rm iness}(s)$, 
$\sharp({\rm non\mbox{-}ori}(s))$ is the number of pole curves of $s$ whose regular neighborhoods are M{\" o}bius bands, 
and $\iota (l)$ is the index of $l$, and we set  $d_0=1$.  

The {\it normalized double bracket polynomial} $R_{(F, D_F)}$ is defined by 
$$ R_{(F, D_F)} = ((-A)^{-3\omega(D)}\dl (F, D_F)\dr 
\quad \in \Z[A, A^{-1}, M, d_1, d_2, \dots],$$ 
where $\omega(D)$ is the writhe of $D$.  
}\end{definition}

For a twisted link diagram $D$, let $(F, D_F)$ be a link diagram realization.  
Then the double bracket polynomial $\dl (F, D_F) \dr$ 
and the normalized double bracket polynomial $ R_{(F, D_F)}$ 
of $(F, D_F)$ are exactly equal to 
the  double bracket polynomial $\dl D \dr$ 
and the normalized double bracket polynomial $ R_D$ 
of $D$ defined in \cite{rkamF}, respectively.

\begin{thm}[\cite{rkamF}]\label{thm4}
The normalized double bracket polynomial $ R_D$ is an invariant of a twisted link.
\end{thm}

For a twisted link $L$, we denote by $R_L$ the polynomial $R_D$, where $D$ is a twisted link diagram representing $L$. 

By Theorem~\ref{thm4}, the following is a corollary to Theorem~\ref{thm1}.  

\begin{cor}\label{cor2}
Let $L$ be a twisted link. Suppose that there is at least one term which includes $d_1^{k_1}\dots d_m^{k_m}$ in the  polynomial invariant $R_L$ of $L$.  
For any link diagram realization $(F, D_F)$  of a link diagram $D$ of $L$, there is a
surface pole state $(F, s)$ of $(F, D_F)$ such that there is a pole curve of $s$ which is not a separating curve in $F$.
\end{cor}

\begin{pf}
By Theorem~\ref{thm4}, for any  link diagram realization $(F, D_F)$ of a twisted link diagram $D$ representing $L$, 
the double bracket polynomial $\dl (F, D_F) \dr$ of $(F, D_F)$ has at least one term which includes $d_1^{k_1}\dots d_m^{k_m}$. 
Thus there is a state $(F, s)$ such that  $\prod_{l \in s}d_{\iota(l)} \neq 1$.  Therefore there is a pole curve $l$ in $s$ with 
$\iota(l) >0$.  By Theorem~\ref{thm1}, the pole curve $l$ is not a separating curve in $F$.  
\qed
\end{pf}

\vspace{.5cm}
{\bf Acknowledgements }

The author  would like to thank Seiichi Kamada for his useful suggestion.

\end{document}